\documentclass[10pt]{article}

\usepackage{mathrsfs}
\usepackage{amsfonts}
\usepackage{amsfonts}
\usepackage{graphicx,latexsym,bm,amsmath,amssymb,verbatim,multicol,lscape}

\date{}

\begin{document}
\title{{\bf The number of rational points on a family of varieties over finite fields}
\author{Shuangnian Hu and Shaofang Hong*\\
{\it Mathematical College, Sichuan University, Chengdu 610064, P.R. China} \\
E-mails: sfhong@scu.edu.cn, s-f.hong@tom.com,\\
 hongsf02@yahoo.com (S. Hong); hushuangnian@163.com (S. Hu)}
\thanks{Hong is the corresponding author and was supported
partially by National Science Foundation of China Grant \#11371260.}}
\maketitle  
{\bf Abstract.}
Let $\mathbb{F}_q$ stand for the finite field of odd characteristic $p$ with $q$
elements ($q=p^{n},n\in \mathbb{N} $) and  $\mathbb{F}_q^*$
denote the set of all the nonzero elements of  $\mathbb{F}_{q}$. Let $m$ and $t$
be positive integers. In this paper, by using the Smith normal form of the exponent matrix,
we obtain a formula for the number of rational points on the
variety defined by the following system of equations over $\mathbb{F}_{q}$:
$$
\sum\limits_{j=0}^{t-1}\sum\limits_{i=1}^{r_{j+1}-r_j}
a_{k,r_j+i}x_1^{e^{(k)}_{r_j+i,1}}...x_{n_{j+1}}^{e^{(k)}_{r_j+i,n_{j+1}}}=b_k, \ k=1,...,m.
$$
where the integers $t>0$, $r_0=0<r_1<r_2<...<r_t$, $1\le n_1<n_2<...<n_t$,
$0\leq j\leq t-1$, $b_k\in \mathbb{F}_{q}$, $a_{k,i}\in \mathbb{F}_{q}^{*}$,
$(k=1,...,m, i=1,...,r_t)$, and the exponent of each variable is a positive integer.
Furthermore, under some natural conditions, we arrive at an explicit formula
for the number of the above variety. It extends the results obtained previously by
Wolfmann, Sun, Wang, Song, Chen, Hong, Hu and Zhao et al. Our result also answers
completely an open problem raised by Song and Chen.  

{\it Keywords and phrases:} Finite field, hypersurface, rational point,
exponent matrix, Smith normal form.

\textit{AMS  Subject Classification:} 11T06, 11T71

\section{Introduction and statement of main result}
Let $\mathbb{F}_{q}$ be the finite field of $q$
elements with odd characteristic $p$ ($q=p^{n},n\in \mathbb{N} $
(the set of positive integers)) and $\mathbb{F}_{q}^*$
denote the nonzero elements of  $\mathbb{F}_{q}$.
By $f_i(x_1, ..., x_n)$ $(i=1,...,m)$ we denote
some polynomials with $n$ variables over $\mathbb{F}_{q}$
and $V$ stands for the following algebraic variety
over $\mathbb{F}_{q}$:

$$\left\{
\begin{aligned}
f_1(x_1, ..., x_n)&=0,\\
...&...\\
f_m(x_1, ..., x_n)&=0.
\end{aligned}
\right.
$$
Let $N_q(V)$ denote the number of $\mathbb{F}_q$-rational
points on the algebraic variety $V$ in $\mathbb{F}^n_q$.
That is, $\mathbb{N}_q(V)=\#\{(x_1,...,x_n)
\in\mathbb{F}_{q}^n|f_i(x_1, ..., x_n)=0,\ i=1,...,m\}.$
Especially, we use $N_q(f)$ to denote $N_q(V)$ if $m=1$.

Studying the exact value of $N_q(V)$ is one of the main topics in finite fields.
The degrees $\text{deg}(f_i)$ play an important role in the estimate of $N_q(V)$.
Let $\lceil x\rceil$ denote the least integer
$\geq x$ and $\text{ord}_q$ denote the  additive valuation such that $\text{ord}_qq=1$.
In 1964, Ax \cite{[A]} generalized the Chevalley-Warning theorem by showing that
$$
\text{ord}_q N_q(V)\geq \Big\lceil\frac{n-\sum_{i=1}^{m}
\text{deg}f_i}{\sum_{i=1}^{m}\text{deg}f_i}\Big\rceil.
$$
Later, further works were done by Katz \cite{[K]}, Adolphson-Sperber \cite{[AS1]}
-\cite{[AS2]}, Moreno-Moreno \cite{[MM]} and Wan \cite{[W1]}-\cite{[W3]}.

It is difficult to give an explicit formula for $N_q(V)$ in general.
Finding explicit formula for $N_q(f)$ under certain conditions
has attracted many authors for many years (see, for instance, \cite{[HV]} et al).
It is well known that there exists an explicit formula for $N_q(f)$ with
$\text{deg}(f) \leq2$ in $\mathbb{F}_{q}$ (see, for example, \cite{[IR]}
and \cite{[L]}). One first considered the diagonal hypersurface:
$$
a_1x_1^{e_1} + ... + a_nx_n^{e_n}-b=0, \ 1\leq i\leq n,
\ a_i\in \mathbb{F}_{q}^{*}, \ b\in \mathbb{F}_{q},\  e_i > 0,\eqno(1.1)
$$
and much work has been done to seek for the number of rational points of
the hypersurface (1.1), see, for instance, \cite{[S1]} and \cite{[Wo1]}-\cite{[Wo2]}.
Carlitz \cite{[Ca1]}, Cohen \cite{[Co]} and Hodges \cite{[Hod]}
counted the rational points on the following $k$-linear hypersurface
$$
a_1x_{11}...x_{1k} + a_2x_{21}...x_{2k} + ... + a_nx_{n1}...x_{nk}-b=0,\eqno(1.2)
$$
with $a_i\in \mathbb{F}_{q}^{*}$, $b\in \mathbb{F}_{q}$. Cao \cite{[C]},
Cao and Sun \cite{[CaSu1]} \cite{[CaSu2]} studied the rational points
on the following more general diagonal hypersurface
$$
a_1x_{11}^{e_{11}}...x_{1n_1}^{e_{1n_1}} +a_2x_{21}^{e_{21}}...x_{2n_2}
^{e_{2n_2}}+ ... + a_rx_{r1}^{e_{r1}}...x_{rn_r}^{e_{rn_r}}=0,
$$
with $1\leq i\leq r$, $1\leq j\leq n_i$, $e_{ij}\in \mathbb{N}$,
$a_i\in \mathbb{F}_{q}^{*}$. Clearly, this extends (1.1) and (1.2)
when $b=0$. On the other hand, 
Pan, Zhao and Cao \cite{[PC]} considered the following hypersurface
$(a_1x_1^{m_1}+...+a_nx_n^{m_n})^\lambda -bx_1^{k_1}...x_n^{k_n}=0$
which extended the results of Carlitz in \cite{[Ca2]} and \cite{[Ca3]}.

If $f=a_1x_1^{e_{11}}...x_n^{e_{1n}} + ... + a_sx_1^{e_{s1}}...x_n^{e_{sn}}-b$
with $e_{ij}>0$ and $a_i\in \mathbb{F}_{q}^{*}$ for $1\leq i\leq s$ and $1\leq j\leq n$,
$b\in \mathbb{F}_{q}$, then a formula for $N_q(f)$ was given by Sun \cite{[S2]}.
Moreover, if $s=n$ and $\text{gcd}(\det(e_{ij}),q-1)=1$
($\det(e_{ij})$ represents the determinant of the $n\times n$ matrix
$\big(e_{ij}\big)$, then Sun \cite{[S2]} gave the explicit formula
for the number of rational points as follows:

\begin{align*}
N(f)={\left\{\begin{array}{rl}q^n-(q-1)^n+
\frac{(q-1)^n+(-1)^n(q-1)}{q},& \ \text{if} \ b=0,\\
\frac{(q-1)^n-(-1)^n}{q},  &  \ \text{if} \ b\neq0.
\end{array}\right.}
\end{align*}
Wang and Sun \cite{[WS]} gave the formula for the
number of rational points of the following hypersurface
$$
a_1x_1^{e_{11}} +a_2x_1^{e_{11}} x_2^{e_{22}}+ ... +
a_nx_1^{e_{n1}}x_2^{e_{n2}}...x_n^{e_{nn}}-b=0,
$$
with  $e_{ij}\geq0$, $a_i\in \mathbb{F}_{q}^{*}$, $b\in \mathbb{F}_{q}$.
In 2005, Wang and Sun \cite{[WS051]} extended  the results of
\cite{[S2]} and \cite{[WS]}. Recently, Hu, Hong and Zhao \cite{[HHZ]}
generalized Wang and Sun's results. In fact, they used the
Smith normal form to present a formula for $N_q(f)$ with $f$
being given by:
\begin{equation*}
f=\sum\limits_{j=0}^{t-1}\sum\limits_{i=1}^{r_{j+1}-r_j}
a_{r_j+i}x_1^{e_{r_j+i,1}} ...x_{n_{j+1}}^{e_{r_j+i,n_{j+1}}}-b, \eqno(1.3)
\end{equation*}
where the integers $t>0$, $r_0=0<r_1<r_2<...<r_t$, $1\le n_1<n_2<...<n_t$,
$b\in \mathbb{F}_{q}$, $a_i\in \mathbb{F}_{q}^{*}$ $(1\leq i\leq r_t)$
and the exponents $e_{ij}$ of each variable are positive integers.

On the other hand, Yang \cite{[Y]} followed Sun's method and gave a formula for
the rational points $N_q(V)$ on the following variety $V$ over $\mathbb{F}_{q}$:
$$\left\{
\begin{aligned}
a_{11}x_1^{e_{11}^{(1)}}...x_n^{e_{1n}^{(1)}} + ... +
a_{1s}x_1^{e_{s1}^{(1)}}...x_n^{e_{sn}^{(1)}}-b_1&=0,\\
...&...\\
a_{m1}x_1^{e_{11}^{(m)}}...x_n^{e_{1n}^{(m)}} + ...
+ a_{ms}x_1^{e_{s1}^{(m)}}...x_n^{e_{sn}^{(m)}}-b_m&=0.
\end{aligned}
\right.
$$
Very recently, Song and Chen \cite{[SC]} continued to make use of
Sun's method and obtained a formula for $N_q(V)$ with $V$ being
the variety over $\mathbb{F}_{q}$ defined by:
$$\left\{
\begin{aligned}
\sum_{j=1}^{s_1}a_{1j}x_1^{e_{j1}^{(1)}}...x_{n_1}^{e_{j,{n_1}}^{(1)}}
+\sum_{j=s_1+1}^{s_2}a_{1j}x_1^{e_{j1}^{(1)}}...x_{n_2}^{e_{j,{n_2}}^{(1)}}-b_1&=0,\\
...&...\\
\sum_{j=1}^{s_1}a_{mj}x_1^{e_{j1}^{(m)}}...x_{n_1}^{e_{j,{n_1}}^{(m)}}
+\sum_{j=s_1+1}^{s_2}a_{mj}x_1^{e_{j1}^{(m)}}...x_{n_2}^{e_{j,{n_2}}^{(m)}}-b_m&=0.
\end{aligned}
\right.
$$
Meanwhile, they proposed an open problem. To state this question,
we need to introduce some notation. In what follows, we always let
$t, m, r_1, ..., r_t, n_1, ..., n_t$ be positive integers such that
$r_1<...<r_t$, $1\le n_1<...<n_t$ and $r_0=0$. For any integers $i, j$ and $k$
with $1\le i\le r_t$, $0\leq j\leq t-1$ and $1\le k\le m$, let $e^{(k)}_{ij}>0$
be integers, $b_k\in \mathbb{F}_{q}$ and $a_{ki}\in \mathbb{F}_{q}^{*}$, and let
$f_k(\mathrm{x}):=f_k(x_1,...,x_{n_t})\in \mathbb{F}_q[x_1,...,x_{n_t}]$
be defined by
$$
f_k(\mathrm{x}):=f_k(x_1,...,x_{n_t})=\sum\limits_{i=1}^{r_t}a_{ki}
\mathrm{x}^{E^{(k)}_i}-b_k, \eqno(1.3)
$$
with $E^{(k)}_i$ being the vectors of non-negative integer components of
dimension $n_t$ defined as follows:
$$\left\{
\begin{aligned}
E^{(k)}_1&=(e^{(k)}_{11},...,e^{(k)}_{1n_1},0,...,0),\  \
\mathrm{x}^{E^{(k)}_1}=x_1^{e^{(k)}_{11}}...x_{n_1}^{e^{(k)}_{1n_1}}, \\
......&...... \\
E^{(k)}_{r_1}&=(e^{(k)}_{r_1,1},...,e^{(k)}_{r_1,n_1},0,...,0),\  \
\mathrm{x}^{E^{(k)}_{r_1}}=x_1^{e^{(k)}_{r_1,1}}...x_{n_1}^{e^{(k)}_{r_1,n_1}}, \\
E^{(k)}_{r_1+1}&=(e^{(k)}_{r_1+1,1},...,e^{(k)}_{r_1+1,n_2},0,...,0),\  \
\mathrm{x}^{E^{(k)}_{r_1+1}}=x_1^{e^{(k)}_{r_1+1,1}}...x_{n_2}^{e^{(k)}_{r_1+1,n_2}},\\
......&...... \\
E^{(k)}_{r_2}&=(e^{(k)}_{r_2,1},...,e^{(k)}_{r_2,n_2},0,...,0), \  \
\mathrm{x}^{E^{(k)}_{r_2}}=x_1^{e^{(k)}_{r_2,1}}...x_{n_2}^{e^{(k)}_{r_2,n_2}},\\
......&...... \\
E^{(k)}_{r_{t-1}+1}&=(e^{(k)}_{r_{t-1}+1,1},...,e^{(k)}_{r_{t-1}+1,n_t}), \ \ \mathrm{x}^{E^{(k)}_{r_{t-1}+1}}=x_1^{e^{(k)}_{r_{t-1}+1,1}}...x_{n_t}^{e^{(k)}_{r_{t-1}+1,n_t}},\\
......&...... \\
E^{(k)}_{r_t}&=(e^{(k)}_{r_t,1},...,e^{(k)}_{r_t,n_t}), \  \
\mathrm{x}^{E^{(k)}_{r_t}}=x_1^{e^{(k)}_{r_t,1}}...x_{n_t}^{e^{(k)}_{r_t,n_t}}.\\
\end{aligned}
\right.
$$
The following interesting question was raised in \cite{[SC]}.\\

{\bf Problem 1.1.} \cite{[SC]} Let $f_1(\mathrm{x}), ..., f_m(\mathrm{x})$
be given as in (1.3). What is the formula for the number of rational points
on the variety defined by the following system of equations over $\mathbb{F}_{q}$:
$$\left\{
\begin{aligned}
f_1(\mathrm{x})&=0,\\
...&...\\
f_m(\mathrm{\mathrm{x}})&=0?
\end{aligned}
\right.\eqno(1.4)
$$

When $m=1$, this question has been answered by Hu, Hong and
Zhao \cite{[HHZ]}. When $t=2$, this question was answered
by Song and Chen \cite{[SC]}. However, if $m\ge 2$ and
$t\ge 3$, then this problem has not been solved yet
and is still kept open so far. Note that a more general
question was proposed in \cite{[HHZ]} and a partial answer
to this general question was given in \cite{[HH]}.

In this paper, our main goal is to investigate Problem 1.1.
We will follow and develop the method of \cite{[HHZ]} to
study Problem 1.1. To state the main result, we first
need to introduce some related concept and notation.
For $1\leq k \leq m$, the exponent matrix of $f_k(\mathrm{x})$,
denoted $E_{f_k}$, is defined to be the $r_t\times n_t$ matrix

\begin{equation*}
E_{f_k}:=
\left( \begin{array}{*{15}c}
E^{(k)}_1\\
\vdots\\
E^{(k)}_{r_t}\\
\end{array}\right)_{r_t\times n_t}.
\end{equation*}
For any integer $l$ with $1\leq l\leq t$, let $E^{(l)}_{f_k}$
($1\le k\le m$) be the $r_l\times n_l$ submatrix of $E_{f_k}$
consisting of the first $r_l$ rows and the first $n_l$ columns
of $E_{f_k}$. Furthermore, we define

\begin{equation*}
E^{(l)}:=
\left( \begin{array}{*{15}c}
E^{(l)}_{f_1}\\
\vdots\\
E^{(l)}_{f_m}\\
\end{array}\right)_{mr_l\times n_l}.
\end{equation*}
Then the famous Smith normal form (see \cite{[S]}, \cite{[Hu]}
or Section 2 below) guarantees the existence of unimodular matrices
$U^{(l)}$ of order $mr_l$ and $V^{(l)}$ of order $n_l$ such that
\begin{equation*}
U^{(l)}E^{(l)}V^{(l)}=
\left( \begin{array}{*{15}c}
D^{(l)}&0\\
0&0
\end{array}\right),
\end{equation*}
where $D^{(l)}:={\rm diag}(d^{(l)}_{1}, ..., d^{(l)}_{s_l})$
with all the diagonal elements $d^{(l)}_{1}, ...,d^{(l)}_{s_l}$
being positive integers such that $d^{(l)}_{1}|...|d^{(l)}_{s_l}$.
Throughout pick $\alpha \in\mathbb{F}^*_q$ to be a fixed primitive
element of $\mathbb{F}_q$. For any $\beta\in\mathbb{F}^*_q$, there
exists exactly an integer $r\in [1, q-1]$ such that $\beta=\alpha^r$.
Such an integer $r$ is called {\it index} of $\beta$ with
respect to the primitive element $\alpha $ (or called the
{\it logarithm of $\beta$ w.r.t. base $\alpha$}), and is denoted by
$\text{ind}_\alpha(\beta):=r$.

Let $k$ and $l$ be integers with $1\leq k\leq m$ and $1\leq l\leq t$.
For the variety defined by the following system of linear equations
over $\mathbb{F}_q^*$:
$$
{\left\{\begin{array}{rl}
&\sum\limits_{i=1}^{r_l}a_{1i}u_{1i}=b_1,\\
&......\\
&\sum\limits_{i=1}^{r_l}a_{mi}u_{mi}=b_m,
\end{array}\right.}\eqno (1.5)
$$
where $a_{k1}, ..., a_{kr_t}, b_k\in {\mathbb{F}_q^*}$
being given in (1.4), we use $N_l$ to denote the number of rational points
$(u_{11}, ..., u_{1r_l}, ..., u_{m1},..., u_{mr_l})\in ({\mathbb{F}}_{q}^{*})^{mr_l}$
of (1.5) under the following extra conditions:
$$
{\left\{\begin{array}{rl}\text{gcd}(q-1,d_i^{(l)})|h_i'& {\rm for} \ i=1,...,s_l,\\
(q-1)|h_i'&{\rm for} \  i=s_l+1,...,mr_l,
\end{array}\right.}\eqno (1.6)
$$
where
$$(h'_1, ...,h'_{mr_l})^T:=U^{(l)}
(\text{ind}_\alpha(u_{11}),..., \text{ind}_\alpha(u_{1r_l})
,...,\text{ind}_\alpha(u_{m1}),...,
\text{ind}_\alpha(u_{mr_l}))^T.$$
Note that by Lemma 2.5 below, one knows that $N_l$ is independent
of the choice of the primitive element $\alpha$.
In what follows, we let $N^{(n_1,...,n_t)}$ stand for the number of
rational points on the variety defined by (1.4). Now let
$$N_{r_0}^{(n_1,...,n_t)}:=q^{n_t-n_1}(q^{n_1}-(q-1)^{n_1}),\eqno(1.7)$$
$$N_{r_t}^{(n_1,...,n_t)}:=N_t(q-1)^{n_t-s_t}
\prod\limits_{j=1}^{s_t}\gcd(q-1, d_j^{(t)}) \eqno(1.8)$$
and for $1\le l\le t-1$, let

$$N_{r_l}^{(n_1,...,n_t)}:=N_lq^{n_t-n_{l+1}}(q-1)^{n_l-s_l}(q^{n_{l+1}-n_l}-
(q-1)^{n_{l+1}-n_l})\prod\limits_{j=1}^{s_l}\gcd(q-1,d_j^{(l)}). \eqno(1.9)$$
Then we are in a position to state the main result of this paper.\\

{\bf Theorem 1.2.} {\it We have
\begin{align*}
N^{(n_1,...,n_t)}={\left\{\begin{array}{rl}
\sum\limits_{i=0}^{t}N_{r_i}^{(n_1,...,n_t)},& \  if  \ b_1=...=b_m=0,\\
\sum\limits_{i=1}^{t}N_{r_i}^{(n_1,...,n_t)}, & \ otherwise,
\end{array}\right.}
\end{align*}
with $N_{r_i}^{(n_1,...,n_t)} (0\le i\le t)$ being defined as in (1.7) to (1.9).}\\
\\
Evidently, Theorem 1.2 extends the main results of \cite{[HHZ]} and \cite{[SC]},
and answers completely Problem 1.1.

This paper is organized as follows. In Section 2, we recall some useful
known lemmas which will be needed later. Subsequently, in Section 3,
we make use of the results presented in Section 2 to show Theorem 1.2.
Also we supply some interesting corollaries. Finally, in Section 4,
we provide two examples to demonstrate the validity of Theorem 1.2.

Throughout this paper, we use $\gcd(a,m)$ to denote the greatest
common divisor of any positive integers $a$ and $m$.

\section{Preliminary lemmas}
In this section, we present some useful lemmas that are needed in section 3.
We first recall two well known definitions.\\

{\bf Definition 2.1.} \cite{[Hu]} Let $M$ be a square integer matrix. If
the determinant of $M$ is $\pm 1$, then $M$ is called a {\it unimodular matrix}.\\

{\bf Definition 2.2.} \cite{[Hu]}  For given any positive integers $m$ and $n$, let
$P$ and $Q$ be two $m\times n$ integer matrices. Suppose that there are two
modular matrices $U$ of order $m$ and  $V$ of order $n$ such that
$P=UQV$. Then we say that $P$ and $Q$ are {\it equivalent}
and we write $P\thicksim Q$.\\

Clearly, the equivalence has the three properties of being reflexive, symmetric and transitive.\\

{\bf Lemma 2.1.} \cite{[Hu]} \cite{[S]}
{\it Let $P$ be a nonzero $m\times n$ integer matrix.
Then $P$ is equivalent to a block matrix of the following form
\begin{equation*}
\left( \begin{array}{*{15}c}
D&0\\
0&0
\end{array}\right), \eqno(2.1)
\end{equation*}
where $D={\rm diag}(d_1, ..., d_r)$ with all the diagonal elements $d_i$
being positive integers and satisfying that $d_i|d_{i+1}$ $(1 \le i < r)$.
In other words, there are unimodular matrices $U$ of order $m$ and $V$
of order $n$ such that}

$$
UPV=\left( \begin{array}{*{15}c}
{\rm diag}(d_1, ..., d_r)&0\\
0&0
\end{array}\right).
$$

We call the diagonal matrix in (2.1) the {\it Smith normal form} of the matrix $P$.
Usually, one writes the Smith normal form of $P$ as ${\rm SNF}(P)$.
The elements $d_i$ are unique up to multiplication by a unit and are called
the {\it elementary divisors, invariants,} or {\it invariant factors}.

For any system of linear congruences
$$
{\left\{\begin{array}{rl}
\sum_{j=1}^{n}h_{1j}y_j &\equiv b_1\pmod m,\\
............\\
\sum_{j=1}^{n}h_{sj}y_j &\equiv b_s\pmod m,
\end{array}\right.}  \eqno (2.2)
$$
let $Y=(y_1,..., y_n)^T$ be the column of indeterminates $y_1, ..., y_n$,
$B=(b_1,...,b_s)^T$ and $H=(h_{ij})$ be the matrix of its
coefficient. Then one can write (2.2) as
$$HY\equiv B\pmod m. \eqno(2.3)$$
By Lemma 2.1, there are unimodular matrices $U$ of order $s$ and $V$
of order $n$ such that
$$
UHV={\rm SNF}(H)=\left( \begin{array}{*{15}c}
{\rm diag}(d_1, ..., d_r)&0\\
0&0
\end{array}\right).
$$

Consequently, we have the following lemma.\\

{\bf Lemma 2.2.} \cite{[HHZ]}
{\it Let $B'=(b'_1,...,b'_s)^T=UB$. Then the system (2.3) of linear congruences is solvable
if and only if $\gcd(m, d_i)|b'_i$ for all integers $i$ with $1\le i\le r$ and $m|b'_i$
for all integers $i$ with $r+1\le i\le s$. Further, the number of solutions of (2.3)
is equal to $m^{n-r}\prod_{i=1}^r\gcd(m, d_i).$}\\

The following result is due to Sun \cite{[S2]}.\\

{\bf Lemma 2.3.} \cite{[S2]}
{\it Let $c_1,...,c_k\in \mathbb{F}_q^{*} $ and  $c \in\mathbb{F}_{q}$, and
let $N(c)$ denote the number of rational points $(x_1, ..., x_k)\in
(\mathbb{F}_{q}^{*})^k$ on the hypersurface $c_1x_1+...+c_kx_k=c.$ Then}
\begin{align*}
N(c)={\left\{\begin{array}{rl}\frac{(q-1)^k+(-1)^k(q-1)}{q},& \ {\it if} \ c=0,\\
\frac{(q-1)^k-(-1)^k}{q},  & {\it otherwise}.
\end{array}\right.}
\end{align*}\\

{\bf Lemma 2.4.}
{\it Let $c_{ij}\in \mathbb{F}_q^{*}$ for all integers $i$ and $j$ with
$1\le i\le m$ and $1\le j\le k$ and  $c_1, ..., c_m\in\mathbb{F}_{q}$.
Let $N(c_1, ..., c_m)$ denote the number of rational points
$(x_{11}, ..., x_{1k}, ..., x_{m1}, ..., x_{mk})\in (\mathbb{F}_{q}^{*})^{mk}$
on the following variety
$$
{\left\{\begin{array}{rl}
&c_{11}x_{11}+...+c_{1k}x_{1k}=c_1,\\
&......\\
&c_{m1}x_{m1}+...+c_{mk}x_{mk}=c_m.
\end{array}\right.} \eqno(2.4)
$$
Then
$$
N(c_1, ..., c_m)=\frac{(q-1)^r}{q^m}((q-1)^{k-1}+(-1)^k)^r((q-1)^k-(-1)^k)^{m-r},
$$
where $r:=\#\{1\le i\le m| c_i=0\}$.}

{\it Proof.} For $1\le i\le m$, let $N(c_i)$ denote the number of
rational points $(x_{i1}, ..., x_{ik})\in (\mathbb{F}_{q}^{*})^{k}$
on the hypersurface $c_{i1}x_{i1}+...+c_{ik}x_{ik}=c_i$. Since
for any rational points
$(x_{11}, ..., x_{1k}, ..., x_{m1}, ..., x_{mk})\in (\mathbb{F}_{q}^{*})^{mk}$
on the variety (2.4), the involved variables are different
from equation to equation, one has
$$N(c_1, ..., c_m)=\prod_{i=1}^mN(c_i).$$
So Lemma 2.3 applied to $N(c_i)$ gives us the required result.
Hence Lemma 2.4 is proved. \hfill$\Box$\\

{\bf Definition 2.3.} Let $k$ be a positive integer.
We say that the column vector $\left( \begin{array}{*{15}c}
c_1\\
\vdots\\
c_k\\
\end{array}\right)\in  \mathbb{Z}^k$ of dimension $k$ {\it divides}
the column vector $\left( \begin{array}{*{15}c}
d_1\\
\vdots\\
d_k\\
\end{array}\right)\in \mathbb{Z}^k$ of dimension $k$ if $c_i$ divides $d_i$
for all integers $i$ with $1\le i\le k$. The divisibility between two
row integer vectors can be defined similarly. \\

{\bf Lemma 2.5.}
{\it Let $\alpha $, $\beta$ be two primitive elements of $\mathbb{F}^*_q$
and $k$ be a positive integer. Let $u_1, ..., u_k$ be nonzero elements of $\mathbb{F}_q$
and each of  $d_1, ..., d_k$ divides $q-1$. Let $U$ be a unimodular matrix.
Then $(d_1, ..., d_k)^T$ divides $U({\rm ind}_{\alpha}u_1, ..., {\rm ind}_{\alpha}u_k)^T$
if and only if
$(d_1, ..., d_k)^T$ divides $U({\rm ind}_{\beta}u_1, ..., {\rm ind}_{\beta}u_k)^T$.
}

{\it Proof.} First of all, for any integer $i$ with $1\le i\le k$, one has
$$
{\rm ind}_{\alpha}u_i\equiv {\rm ind}_{\alpha}\beta {\rm ind}_{\beta}u_i \pmod {q-1}.
$$
It then follows that
$$
U\left( \begin{array}{*{15}c}
{\rm ind}_{\alpha}u_1\\
\vdots\\
{\rm ind}_{\alpha}u_k\\
\end{array}\right)
\equiv {\rm ind}_{\alpha}\beta \cdot U\left( \begin{array}{*{15}c}
{\rm ind}_{\beta}u_1\\
\vdots\\
{\rm ind}_{\beta}u_k\\
\end{array}\right)  \pmod {q-1}.
\eqno(2.5)
$$

On the one hand, since all of $d_1, ..., d_k$ divide $q-1$, by (2.5) one knows that
$(d_1, ..., d_k)^T$ divides $U({\rm ind}_{\alpha}u_1, ..., {\rm ind}_{\alpha}u_k)^T$
if and only if $(d_1, ..., d_k)^T$ divides ${\rm ind}_{\alpha}\beta \cdot
U({\rm ind}_{\beta}u_1, ..., {\rm ind}_{\beta}u_k)^T$.

On the other hand, since $\alpha $ and $\beta $ are primitive elements,
${\rm ind}_{\alpha}\beta$ is coprime to $q-1$. Hence ${\rm ind}_{\alpha}\beta$
is coprime to each of $d_1, ..., d_k$. Then one can derive
that $(d_1, ..., d_k)^T$ divides ${\rm ind}_{\alpha}\beta \cdot U({\rm ind}_{\beta}u_1,
..., {\rm ind}_{\beta}u_k)^T$ if and only if $(d_1, ..., d_k)^T$ divides
$U({\rm ind}_{\beta}u_1, ..., {\rm ind}_{\beta}u_k)^T$.

Finally, the desired result follows immediately. So Lemma 2.5 is proved. \hfill$\Box$\\

{\bf Remark 2.1.} If $N_l(\alpha )$ stands for the number of rational points
$(u_{11}, ..., u_{1r_l},\\ ..., u_{m1},..., u_{mr_l})\in ({\mathbb{F}}_{q}^{*})^{mr_l}$
of (1.5) under the extra conditions (1.6) with respect to the primitive element $\alpha $,
then by Lemma 2.5 we have that $N_l(\beta)=N_l(\gamma)$ for any primitive elements $\beta$
and $\gamma$. So we can use $N_l$ to denote the number of rational points
$(u_{11}, ..., u_{1r_l}, ..., u_{m1},..., u_{mr_l})\in ({\mathbb{F}}_{q}^{*})^{mr_l}$
of (1.5) under the extra conditions (1.6) with respect to any given primitive element $\alpha $.

\section{Proof of Theorem 1.2}

In this section, we give the proof of Theorem 1.2. First, we present
some notation and two lemmas. For any given $(u_{11},...,u_{1r_t},...,u_{m1},
...,u_{m,r_t})\in \mathbb{ F}_{q}^{mr_t}$, we use
$$N(\mathrm{x}^{E^{(k)}_i}=u_{ki}, k=1,...,m, i=1,..., r_t)$$
to denote the number of rational points $(x_1,...,x_{n_t})
\in\mathbb{ F}_{q}^{n_t}$ on the following algebraic
variety over $\mathbb{ F}_{q}$:
$$
\mathrm{x}^{E^{(k)}_i}=u_{ki}, k=1,...,m, i=1,..., r_t, \eqno(3.1)
$$
Define $T$ to be the set of rational points $(u_{11},...,u_{1r_t},...,
u_{m1},...,u_{m,r_t})\in\mathbb{ F}_{q}^{mr_t}$ on the variety
$\sum_{i=1}^{r_t}a_{ki}u_{ki}=b_k$, $k=1,...,m$. Namely,
$$T:=\Big\{(u_{11},...,u_{1r_t},...,u_{m1},...,u_{m,r_t})\in
\mathbb{ F}_{q}^{mr_t}:\sum_{i=1}^{r_t}a_{ki}u_{ki}=b_k, k=1,...,m\Big\},$$
with $b_k$ and $a_{ki}$ $(1\leqslant k\leqslant m, 1\leqslant i\leqslant r_t)$
being given as in (1.3).

Let $T(0)$ consist of zero vector of dimension $mr_t$.
For any integer $n$ with $1\le n\le r_t$, let $T(n)$ denote the subset
of $T$ in which the vector holds exactly $mn$ nonzero components.
Define
$$M_n^{(n_1,...,n_t)}:=\sum\limits_{(u_{11},...,u_{1r_t},...,u_{m1},...,u_{m,r_t})
\in T(n)}N(\mathrm{x}^{E^{(k)}_i}=u_{ki}, k=1,...,m, i=1, ..., r_t).\eqno(3.2)
$$
Then we have the following lemma.\\

{\bf Lemma 3.1.} {\it One has
$$
N^{(n_1,...,n_t)}=\sum_{n=0}^{r_t}M_n^{(n_1,...,n_t)}.$$
}

{\it Proof.}
First of all, by (1.3) and (3.1) we obtain that
$$
N^{(n_1,...,n_t)}=\sum\limits_{(u_{11},...,u_{1r_t},...,u_{m1},
...,u_{m,r_t})\in \mathbb{ F}_{q}^{mr_t}
\atop a_{k1}u_{k1}+...+a_{{kr_t}}u_{kr_t}=b_k, 1\le k\le m. }
N(\mathrm{x}^{E^{(k)}_i}=u_{ki}, k=1,...,m, i=1,..., r_t).
$$
It follows from the definition of $T$ that
$$
N^{(n_1,...,n_t)}=\sum\limits_{(u_{11},...,u_{1r_t},...,u_{m1},...,u_{m,r_t})
\in T}N(\mathrm{x}^{E^{(k)}_i}=u_{ki}, k=1,...,m, i=1,..., r_t).
\eqno(3.3)
$$

Now let $(u_{11},...,u_{1r_t},...,u_{m1},...,u_{m,r_t})\in T$ and consider
the rational points $(x_1,...,x_{n_t})\in\mathbb{ F}_{q}^{n_t}$ on
the variety (3.1). For any integer $k$ with $1\le k\le m$, one can
easily deduce that for any given integer $i$ with $0<i<r_t$, if $u_{ki}=0$,
then $u_{k,i+1}=0$. Thus if $(u_{11},...,u_{1r_t},...,u_{m1},...,u_{m,r_t})\in T$
is a nonzero vector, then it must hold exactly $mn$ nonzero components
for some integer $n$ with $1\le n\le r_t$. For any integer $n$ with
$1\le n\le r_t$, let $T(n)$ denote the subset of $T$ in which the
vector holds exactly $mn$ nonzero components. So if
$(u_{11},...,u_{1r_t},...,u_{m1},...,u_{m,r_t})\in T(n)$, then
we have that $u_{k1}\neq0, ..., u_{kn}\neq0$ and $u_{k,n+1}=...=u_{k,r_t}=0$
for all integers $k$ with $1\le k\le m$. We claim that if
$(u_{11},...,u_{1r_t},...,u_{m1},...,u_{m,r_t})\in T$ holds exactly $mn$
nonzero components, then $n$ equals one of $r_1, ..., r_t$.
 In the following we will prove the claim. Let $n\in \{1, 2, ..., r_t\}
\setminus \{r_1, ..., r_t\}$. Then one can find an integer $l$
with $1\leq l\leq t$ such that $r_{l-1}<n<r_l$. Assume that $T(n)$ is nonempty.
On the one hand, by the definition of $T(n)$, one has
$u_{k1}\neq0, ..., u_{kn}\neq0$ and $u_{k,n+1}=...=u_{k,r_t}=0$
for all integers $k$ with $1\le k\le m$. On the other hand, from $u_{k,n+1}=0$,
and noting that $u_{kn}=x_1^{e_{n,1}^{(k)}}...x_{n_l}^{e_{n,n_l}^{(k)}}$
and $u_{k,n+1}=x_1^{e_{n+1,1}^{(k)}}...x_{n_l}^{e_{n+1,n_l}^{(k)}}$,
we infer that $u_{kn}=0$. This contradicts to $u_{kn}\neq0$.
So the assumption is not true, and the claim is proved.

Now let $T(0)$ consist of zero vector of dimension $mr_t$.
Then by the claim we have
$$
T=\bigcup_{n=0}^{r_t}T(n).
$$
Thus from (3.3) one derives that
$$
N^{(n_1,...,n_t)}=\sum_{n=0}^{r_t}M_n^{(n_1,...,n_t)}
$$
as required. This ends the proof of Lemma 3.1. \hfill$\Box$\\

Consequently, we compute $M_{r_l}^{(n_1,...,n_t)}$ for integers $l$ with $1\le l\le t$.\\

{\bf Lemma 3.2.} {\it For all integers $l$ with $1\le l\le t$, we have}
$$M_{r_l}^{(n_1,...,n_t)}=N_{r_l}^{(n_1,...,n_t)}.$$

{\it Proof.} Evidently, one has
\begin{align*}
&M_{r_t}^{(n_1,...,n_t)}\\
&=\sum_{(u_{11},...,u_{1r_t},...,u_{m1},...,u_{m,r_t})
\in T(r_t)}N(\mathrm{x}^{E^{(k)}_i}=u_{ki}, k=1,...,m,\ i=1,..., r_t)\\
&=\sum_{\begin{subarray}{I}(u_{11},...,u_{1r_t},...,u_{m1},
...,u_{m,r_t})\in (\mathbb{ F}_{q}^*)^{mr_t}\\
a_{k,1}u_{k1}+...+a_{{k,r_t}}u_{kr_t}=b_k \\
\end{subarray}}
N(\mathrm{x}^{E^{(k)}_i}=u_{ki}, k=1,...,m, i=1,..., r_t). \ (3.4)
\end{align*}

Now we let $l$ be a given integer with
$1\le l\le t-1$. Then
\begin{align*}
&M_{r_l}^{(n_1,...,n_t)}\\
&=\sum_{(u_{11},...,u_{1r_t},...,u_{m1},...,u_{m,r_t})
\in T(r_l)}N(\mathrm{x}^{E^{(k)}_i}=u_{ki}, k=1,...,m,\ i=1,..., r_t)\\
&=\sum_{\begin{subarray}{I}(u_{11},...,u_{1r_t},...,u_{m1},...,u_{m,r_t})\in \mathbb{ F}_{q}^{mr_t}\\
a_{k,1}u_{k1}+...+a_{{k,r_l}}u_{kr_l}=b_k,\ u_{ki}\in \mathbb{ F}_{q}^* \\
u_{kj}=0, k=1,...,m. i=1,...r_l, j=r_l+1,...,r_t \end{subarray}}
N(\mathrm{x}^{E^{(k)}_i}=u_{ki}, k=1,...,m, i=1,..., r_t)\\
&=\sum_{\begin{subarray}{I}(u_{11},...,u_{1r_l},...,u_{m1},...,u_{mr_l})\in (\mathbb{ F}_{q}^*)^{mr_l}\\
a_{k,1}u_{k1}+...+a_{{k,r_l}}u_{kr_l}=b_k,\  k=1,...,m.\end{subarray}}
N\big(\begin{subarray}{I}\mathrm{x}^{E^{(k)}_i}=u_{ki}\ {\rm and}\ \mathrm{x}^{E^{(k)}_j}=0, \\
\ {\rm for}\ 1\leq k\leq m, 1\leq i\leq r_l \ {\rm and}\ r_l+1\leq j
\leq r_t\end{subarray}\big). \ \ \ \ \ \ \ \ \ \ (3.5)
\end{align*}
But the definition of $\mathrm{x}^{E^{(k)}_j}$ tells us that the fact that
$\mathrm{x}^{E^{(k)}_j}=0$ for $r_l+1\le j\le r_t$
is reduced to saying that $\mathrm{x}^{E^{(1)}_{r_l+1}}=0$. It then follows from (3.5) that
\begin{align*}
&M_{r_l}^{(n_1,...,n_t)}\\
&=\sum_{\begin{subarray}{I}(u_{11},...,u_{1r_l},...,u_{m1},...,u_{mr_l})\in (\mathbb{ F}_{q}^*)^{mr_l}\\
a_{k,1}u_{k1}+...+a_{{k,r_l}}u_{kr_l}=b_k,\  k=1,...,m.\end{subarray}}
N\big(\begin{subarray}{I}\mathrm{x}^{E^{(k)}_i}=u_{ki}
\ {\rm for}\ 1\leq k\leq m,\  1\leq i\leq r_l
\ {\rm and}\ \mathrm{x}^{E^{(1)}_{r_l+1}}=0\end{subarray}\big). \ \ (3.6)
\end{align*}

It is easy to see that $\mathrm{x}^{E^{(1)}_{r_l+1}}=0$ is equivalent to
$x_1...x_{n_l}x_{n_l+1}...x_{n_{l+1}}=0$. Since $u_{kr_l}\ne 0$ and
$u_{kr_l}= x_1^{e^{(k)}_{r_l,1}}...x_{n_l}^{e^{(k)}_{r_l,n_l}}$,
one has $x_1...x_{n_l}\ne 0$. So $\mathrm{x}^{E^{(1)}_{r_l+1}}=0$
is equivalent to $x_{n_l+1}...x_{n_{l+1}}=0$.
Then by (3.6), one gets that for $1\le l\le t-1$,
\begin{align*}
&M_{r_l}^{(n_1,...,n_t)}\\
&=\sum_{\begin{subarray}{I}(u_{11},...,u_{1r_l},...,u_{m1},...,u_{mr_l})\in (\mathbb{ F}_{q}^*)^{mr_l}\\
a_{k,1}u_{k1}+...+a_{{k,r_l}}u_{kr_l}=b_k,\  k=1,...,m.\end{subarray}}
N\big(\begin{subarray}{I}\mathrm{x}^{E^{(k)}_i}=u_{ki},
1\leq k\leq m, 1\leq i\leq r_l\ {\rm and} \ x_{n_l+1}...x_{n_{l+1}}=0\end{subarray}\big). \ (3.7)
\end{align*}

For any given $(u_{11},...,u_{1r_l},...,u_{m1},...,u_{mr_l})\in (\mathbb{ F}_{q}^*)^{mr_l}$
with $\sum_{i=1}^{r_l}a_{ki}u_{ki}=b_k (k=1, ..., m)$, one has
\begin{align*}
&N(\mathrm{x}^{E^{(k)}_i}=u_{ki}, 1\leq k\leq m, 1\leq i\leq r_l\ {\rm and} \ x_{n_l+1}...x_{n_{l+1}}=0)=\\
&\#\{(x_1,...,x_{n_t})\in (\mathbb{F}_q)^{n_t}: \mathrm{x}^{E^{(k)}_i}=u_{ki},
1\leq k\leq m, 1\leq i\leq r_l\ {\rm and} \ x_{n_l+1}...x_{n_{l+1}}=0\}.
\end{align*}
Since each of the components $x_{n_{l+1}+1},..., x_{n_t}$ can run over
the whole finite field $\mathbb{F}_q$ independently, it then follows that
\begin{align*}
&N(\mathrm{x}^{E^{(k)}_i}=u_{ki}, 1\leq k\leq m, 1\leq i\leq r_l\ {\rm and} \ x_{n_l+1}...x_{n_{l+1}}=0)\\
&=q^{n_t-n_{l+1}}\times \#\{(x_1,...,x_{n_{l+1}})\in (\mathbb{F}_q)^{n_{l+1}}:\\
&\mathrm{x}^{E^{(k)}_i}=u_{ki},
1\leq k\leq m, 1\leq i\leq r_l\ {\rm and} \ x_{n_l+1}...x_{n_{l+1}}=0\}. \ \ \ \ \ \ (3.8)
\end{align*}
Notice that the choice of $(x_1,...,x_{n_{l}})\in (\mathbb{F}_q^*)^{n_{l}}$
satisfying that $\mathrm{x}^{E^{(k)}_i}=u_{ki}$ ($k=1,...,m$, $i=1,...,r_l$)
is independent of the choice of
$(x_{n_l+1},...,x_{n_{l+1}})\in (\mathbb{F}_q)^{n_{l+1}-n_l}$
satisfying that $x_{n_l+1}...x_{n_{l+1}}=0$. We then derive that
\begin{align*}
&\#\{(x_1,...,x_{n_{l+1}})\in (\mathbb{F}_q)^{n_{l+1}}: \mathrm{x}^{E^{(k)}_i}=u_{ki},
1\leq k\leq m, 1\leq i\leq r_l\ {\rm and} \ x_{n_l+1}...x_{n_{l+1}}=0\}\\
&=\#\{(x_1,...,x_{n_{l}})\in (\mathbb{F}_q^*)^{n_{l}}: \mathrm{x}^{E^{(k)}_i}=u_{ki},
1\leq k\leq m, 1\leq i\leq r_l\}\times \\
&\#\{(x_{n_l+1},...,x_{n_{l+1}})\in (\mathbb{F}_q)^{n_{l+1}-n_l}:
x_{n_l+1}...x_{n_{l+1}}=0\}.
\ \ \ \ \ \ \ \ \ \ \ \ \ \ \ \ \ \ \ \ \ \ \ \ \ \ \  (3.9)
\end{align*}

On the other hand, we can easily compute that
$$
\#\{(x_{n_l+1},...,x_{n_{l+1}})\in (\mathbb{F}_q)^{n_{l+1}-n_l}: x_{n_l+1}...x_{n_{l+1}}=0\}=
$$
$$
\sum_{i=1}^{n_{l+1}-n_l}{n_{l+1}-n_l \choose i}(q-1)^{n_{l+1}-n_l-i}=q^{n_{l+1}-n_l}-(q-1)^{n_{l+1}-n_l}. \eqno(3.10)
$$
So by (3.8) to (3.10), one obtains that
\begin{align*}
&N(\mathrm{x}^{E^{(k)}_i}=u_{ki}, 1\leq k\leq m, 1\leq i\leq r_l\ {\rm and} \ x_{n_l+1}...x_{n_{l+1}}=0)\\
&=q^{n_t-n_{l+1}}(q^{n_{l+1}-n_l}-(q-1)^{n_{l+1}-n_l})\\
&\times \#\{(x_1,...,x_{n_{l}})\in (\mathbb{F}_q^*)^{n_{l}}: \mathrm{x}^{E^{(k)}_i}=u_{ki},
1\leq k\leq m, 1\leq i\leq r_l\}\\
&=q^{n_t-n_{l+1}}(q^{n_{l+1}-n_l}-(q-1)^{n_{l+1}-n_l}) N(\mathrm{x}^{E^{(k)}_i}=u_{ki},
1\leq k\leq m, 1\leq i\leq r_l). \   (3.11)
\end{align*}
Then by (3.7) together with (3.11), we have
\begin{align*}
&M_{r_l}^{(n_1,...,n_t)}=q^{n_t-n_{l+1}}(q^{n_{l+1}-n_l}-(q-1)^{n_{l+1}-n_l})\times \\
&\sum_{\begin{subarray}{I}(u_{11},...,u_{1r_l},...,u_{m1},...,u_{mr_l})\in (\mathbb{ F}_{q}^*)^{mr_l}\\
a_{k,1}u_{k1}+...+a_{{k,r_l}}u_{kr_l}=b_k,\  k=1,...,m.\end{subarray}}
N\big(\mathrm{x}^{E^{(k)}_i}=u_{ki},
1\leq k\leq m, 1\leq i\leq r_l\big). \ \ \ \ \   (3.12)
\end{align*}

Now we treat with the sum
$$
\sum_{\begin{subarray}{I}(u_{11},...,u_{1r_l},...,u_{m1},...,u_{mr_l})
\in (\mathbb{ F}_{q}^*)^{mr_l}\\
a_{k,1}u_{k1}+...+a_{{k,r_l}}u_{kr_l}=b_k,\  k=1,...,m.\end{subarray}}
N\big(\mathrm{x}^{E^{(k)}_i}=u_{ki},
1\leq k\leq m, 1\leq i\leq r_l\big),
$$
where $l=1,...,t$. First, for any given
$(u_{11},...,u_{1r_l},...,u_{m1},...,u_{mr_l})\in (\mathbb{ F}_{q}^*)^{mr_l}$
with $\sum_{i=1}^{r_l}a_{ki}u_{ki}=b_k$ $(k=1,...,m)$,
$N\big(\begin{subarray}{I}\mathrm{x}^{E^{(k)}_i}=u_{ki},
1\leq k\leq m, 1\leq i\leq r_l\end{subarray}\big)$ equals the number
of rational points $(x_1,..., x_{n_l})\in (\mathbb{F}^*_q)^{n_l}$
on the following algebraic variety:

$$\left\{
\begin{aligned}
x_1^{e^{(1)}_{11}}...x_{n_1}^{e^{(1)}_{1n_1}}&=u_{11}, \\
............ \\
x_1^{e^{(1)}_{r_l,1}}...x_{n_l}^{e^{(1)}_{r_l,n_l}}&=u_{1, r_l},\\
............ \\
x_1^{e^{(m)}_{11}}...x_{n_1}^{e^{(m)}_{1n_1}}&=u_{m1}, \\
............ \\
x_1^{e^{(m)}_{r_l,1}}...x_{n_l}^{e^{(m)}_{r_l,n_l}}&=u_{m, r_l}.
\end{aligned}
\right.\eqno(3.13)
$$

Since $u_{k1}...u_{kr_l}\neq0$ $(1\leq k\leq m)$,
we infer that the number of the rational points
$(x_1,..., x_{n_l})\in (\mathbb{F}^*_q)^{n_l}$ of (3.13) is equal to
the number of nonnegative integral solutions
$({\rm ind}_{\alpha}(x_1), ...,{\rm ind}_{\alpha}(x_{n_l}))\in \mathbb{N}^{n_l}$
of the following system of congruences
$$\left\{
\begin{aligned}
\sum\limits_{i=1}^{n_l}e^{(1)}_{1i}\text{ind}_\alpha(x_i)
&\equiv\text{ind}_\alpha (u_{11})\pmod{q-1}, \\
......&...... \\
\sum\limits_{i=1}^{n_l}e^{(1)}_{r_l,i}\text{ind}_\alpha(x_i)
&\equiv\text{ind}_\alpha (u_{1,r_l})\pmod{q-1},\\
......&...... \\
\sum\limits_{i=1}^{n_l}e^{(m)}_{1i}\text{ind}_\alpha(x_i)
&\equiv\text{ind}_\alpha (u_{m1})\pmod{q-1}, \\
......&...... \\
\sum\limits_{i=1}^{n_l}e^{(m)}_{r_l,i}\text{ind}_\alpha(x_i)
&\equiv\text{ind}_\alpha (u_{m, r_l})\pmod{q-1}.
\end{aligned}
\right. \eqno(3.14)
$$
But Lemma 2.2 tells us that (3.14) has
solutions $({\rm ind}_{\alpha}(x_1), ..., {\rm ind}_{\alpha}(x_{n_l}))
\in \mathbb{N}^{n_l}$ if and only if the extra conditions (1.6) hold.
Further, Lemma 2.2 gives us the number of
solutions $({\rm ind}_{\alpha}(x_1), ..., {\rm ind}_{\alpha}(x_{n_l}))
\in \mathbb{N}^{n_l}$ of (3.14) which is equal to
$(q-1)^{n_l-s_l}\prod\limits_{i=1}^{s_l}\gcd(q-1,d_i^{(l)})$. Hence
$$
N\big(\begin{subarray}{I}E^{(k)}_i=u_{ki},
1\leq k\leq m, 1\leq i\leq r_l\end{subarray}\big)=(q-1)^{n_l-s_l}
\prod\limits_{i=1}^{s_l}\gcd(q-1, d_i^{(l)}).\eqno(3.15)
$$

Notice that
$$
N_l=
\sum_{\begin{subarray}{I}(u_{11},...,u_{1r_l},...,u_{m1},...,u_{mr_l})
\in (\mathbb{ F}_{q}^*)^{mr_l}\\
a_{k,1}u_{k1}+...+a_{{k,r_l}}u_{kr_l}=b_k\\
  k=1,...,m \ {\rm and} \ (1.6) \ {\rm holds}\end{subarray}}1.
$$
It then follows from (3.15) that
\begin{align*}
&\sum_{\begin{subarray}{I}(u_{11},...,u_{1r_l},...,u_{m1},...,u_{mr_l})\in (\mathbb{ F}_{q}^*)^{mr_l}\\
a_{k,1}u_{k1}+...+a_{{k,r_l}}u_{kr_l}=b_k,\  k=1,...,m.\end{subarray}}
N\big(\begin{subarray}{I}\mathrm{x}^{E^{(k)}_i}=u_{ki},
1\leq k\leq m, 1\leq i\leq r_l\end{subarray}\big) \\
&=\sum_{\begin{subarray}{I}(u_{11},...,u_{1r_l},...,u_{m1},...,u_{mr_l})\in (\mathbb{ F}_{q}^*)^{mr_l}\\
a_{k,1}u_{k1}+...+a_{{k,r_l}}u_{kr_l}=b_k\\
  k=1,...,m \ {\rm and} \ (1.6) \ {\rm holds}\end{subarray}}
(q-1)^{n_l-s_l}\prod\limits_{i=1}^{s_l}\gcd(q-1, d_i^{(l)})\\
&=(q-1)^{n_l-s_l}\prod\limits_{i=1}^{s_l}\gcd(q-1, d_i^{(l)})
\sum_{\begin{subarray}{I}(u_{11},...,u_{1r_l},...,u_{m1},...,u_{mr_l})\in (\mathbb{ F}_{q}^*)^{mr_l}\\
a_{k,1}u_{k1}+...+a_{{k,r_l}}u_{kr_l}=b_k\\
  k=1,...,m \ {\rm and} \ (1.6) \ {\rm holds}\end{subarray}}1\\
&=N_l(q-1)^{n_l-s_l}\prod\limits_{i=1}^{s_l}\gcd(q-1, d_i^{(l)}). \ \ \ \ \ \ \ \ \ \ \ \  \ \ \ \ \ \ \ \
\ \ \ \ \ \ \ \ \ \ \ \  \ \ \ \ \ \ \ \ \ \ \ \ \ \ \ \ \ (3.16)
\end{align*}

For $1\leq l\leq t$, using (3.3), (3.12) and (3.16), we obtain
the desired result $M_{r_l}^{(n_1,...,n_t)}=N_{r_l}^{(n_1,...,n_t)}$.
This concludes the proof of Lemma 3.2. \hfill$\Box$\\

We can now turn our attention to the proof of Theorem 1.2.\\

{\it Proof of Theorem 1.2.}
First we show that the following is true:
$$
M_{0}^{(n_1,...,n_t)}:={\left\{\begin{array}{rl}
N_{r_0}^{(n_1,...,n_t)},&  \  \text{if}  \ b_1=...=b_m=0,\\
0,  & \ \text{otherwise}.
\end{array}\right.} \eqno(3.17)
$$

Let $b_i\neq0$ for some integer $i$ with $1\le i\le m$.
Then the following variety
$$\left\{
\begin{aligned}
&a_{11}u_{11}+...+a_{1, r_t}u_{1, r_t}=b_1 \\
&.............. \\
&a_{m1}u_{m1}+...+a_{m, r_t}u_{m, r_t}=b_m
\end{aligned}
\right.
$$
does not contain the original point $(0, ..., 0)\in \mathbb{ F}_{q}^{mr_t}$.
So $T(0)$ is empty. It then follows from (3.2) that
$M_{0}^{(n_1,...,n_t)}=0.$
That is, the second part of (3.17) is true.

Now let $b_1=...=b_m=0$. Then $T(0)$ consists of zero vector of
dimension $mr_t$. Thus by (3.2), we have
\begin{align*}
M_{0}^{(n_1,...,n_t)}&=N(\mathrm{x}^{E^{(k)}_i}=0, k=1,...,m, i=1,2,..., r_t)\\
&=N(\mathrm{x}^{E^{(1)}_1}=0)\\
&=N(x_1^{e^{(1)}_{11}}...x_{n_1}^{e^{(1)}_{1n_1}}=0)\\
&=q^{n_t-n_1}\sum_{j=1}^{n_1}\binom{n_1}{j}(q-1)^{n_1-j}\\
&=q^{n_t-n_1}(q^{n_1}-(q-1)^{n_1})=N_{r_0}^{(n_1,...,n_t)}
\end{align*}
as required. This completes the proof of (3.17).

Consequently, we show that for all integers $n\in \{1, ..., r_t\}
\setminus \{r_1, ..., r_t\}$, one has
$$M_n^{(n_1,...,n_t)}=0. \eqno(3.18)$$
To prove it, we choose an $n\in \{1, 2, ..., r_t\}
\setminus \{r_1, ..., r_t\}$. Then there exists an integer $l$
with $1\leq l\leq t$ such that $r_{l-1}<n<r_l$. Claim that
$T(n)$ is empty which will be shown in what follows.

Suppose that $T(n)$ is nonempty. On the one hand, the definition of
$T(n)$ gives us that for all integers $k$ with $1\le k\le m$, we have
$u_{k1}\neq0$,...,$u_{kn}\neq0$ and $u_{k,n+1}=...=u_{k,r_t}=0.$
On the other hand, from $u_{k,n+1}=0$, and noting that
$$u_{kn}=x_1^{e^{(k)}_{n,1}}...x_{n_l}^{e^{(k)}_{n,n_l}},
u_{k,n+1}=x_1^{e^{(k)}_{n+1,1}}...x_{n_l}^{e^{(k)}_{n+1,n_l}},$$
we deduce that $u_{kn}=0$. This contradicts with the fact
$u_{kn}\neq0$. So the assumption is not true.
Hence $T(n)=\varnothing$. Thus by (3.2),
$$M_n^{(n_1,...,n_t)}=0.$$
This finishes the proof of (3.18).

Finally, using Lemmas 3.1 and 3.2, (3.17) and (3.18), the desired
result follows immediately. So Theorem 1.2 is proved. \hfill$\Box$
\\

In concluding this section, we present an interesting corollary.\\

{\bf Corollary 3.1.} {\it If ${\rm SNF}(E^{(l)})=(D^{(l)} \ \ 0)$ and
$\gcd(\det D^{(l)}, q-1)=1$ for all integers $l$ with $1\le l\le t$,
then the number $N^{(n_1,...,n_t)}$ of rational points on the variety (1.4) is
given by
\begin{align*}
N^{(n_1,...,n_t)}={\left\{\begin{array}{rl}
\sum\limits_{i=0}^{t}\widehat{N}_{r_i}^{(n_1,...,n_t)},&{\it if} \  b_1=...=b_m=0,\\
\sum\limits_{i=1}^{t}\widehat{N}_{r_i}^{(n_1,...,n_t)}, & \ {\it otherwise},
\end{array}\right.}
\end{align*}
where $\widehat{N}_{r_0}^{(n_1,...,n_t)}:=q^{n_t-n_1}(q^{n_1}-(q-1)^{n_1})$,
$\widehat{N}_{r_t}^{(n_1,...,n_t)}:=\widehat{N}_t(q-1)^{n_t-s_t}$ and
$\widehat{N}_{r_l}^{(n_1,...,n_t)}:=q^{n_t-n_{l+1}}
(q^{n_{l+1}-n_l}-(q-1)^{n_{l+1}-n_l}) \widehat{N}_l(q-1)^{n_l-s_l}$
$(l=1,...,t-1)$, and for all integers $k$ with $1\leqslant k\leqslant t$, one has
$$
\widehat{N}_k:=\frac{(q-1)^r}{q^m}((q-1)^{r_k-1}+(-1)^{r_k})^r((q-1)^{r_k}-(-1)^{r_k})^{m-r},
$$
where $r:=\#\{1\le i\le m| b_i=0\}$.}

{\it Proof.} Since ${\rm SNF}(E^{(l)})=(D^{(l)} \ \ 0)$ and
$\gcd(\det D^{(l)}, q-1)=1$ for all integers $l$ with $1\le l\le t$,
we derive that
$$s_l=mr_l \ {\rm and} \ \gcd(q-1, d_j^{(l)})=1$$
for all integers $j$ with $1\le j\le s_l$. Thus for any solution
$(u_{11},...,u_{1,r_l},...,u_{m1},...,\\u_{m,r_l})
\in ({\mathbb{F}}_{q}^{*})^{mr_l}$
of (1.5), the extra conditions (1.6) are always true.
It then follows from Theorem 1.2 and Lemma 2.4 that
the desired result follows immediately.
This finishes the proof of Corollary 3.1.\hfill$\Box$

\section {Two examples}
In this section, we supply two examples to illustrate the validity of our main result.\\

{\bf Example 4.1.} We use Corollary 3.1 to compute the number $N^{(n_1,n_2)}$ of
rational points on the following variety over $\mathbb{F}_{11}$:
\begin{align*}
{\left\{\begin{array}{rl}
x_1x_2^3x_3^2+x_1^5x_2^7x_3^5x_4^5x_5x_6^2+x_1^5x_2^4x_3^3x_4^2x_5^6x_6^3=4,\\
x_1x_2^5x_3^3+x_1^3x_2^5x_3^6x_4^5x_5^4x_6^7+x_1^3x_2x_3^5x_4^7x_5^3x_6^7=0.
\end{array}\right.}
\end{align*}

Clearly, we have $n_1=3$, $n_2=6$, $r_1=1$, $r_2=3$, $m=2$,
$$E^{(1)}=\left( \begin{array}{*{20}c}
1&3&2\\
1&5&3
\end{array}\right)
\ \text{and} \
E^{(2)}=\left( \begin{array}{*{20}c}
1&3&2&0&0&0\\
5&7&5&5&1&2\\
5&4&3&2&6&3\\
1&5&3&0&0&0\\
3&5&6&5&4&7\\
3&1&5&7&3&7
\end{array}\right).$$

One can easily deduce that the Smith normal forms of
$E^{(1)}$  and $E^{(2)}$ are given as follows:
$${\rm SNF}(E^{(1)})=
\left( \begin{array}{*{20}c}
1&0&0\\
0&1&0
\end{array}\right)
\ \text{and} \
{\rm SNF}(E^{(2)})=\left( \begin{array}{*{20}c}
1&0&0&0&0&0\\
0&1&0&0&0&0\\
0&0&1&0&0&0\\
0&0&0&1&0&0\\
0&0&0&0&1&0\\
0&0&0&0&0&291
\end{array}\right).$$
Thus $s_1=2$ and $s_2=6$. It follows from Corollary 3.1 that
$$N^{(n_1,n_2)}=\widehat{N}_{r_1}^{(n_1,n_2)}+\widehat{N}_{r_2}^{(n_1,n_2)}
=0+\frac{(10^3+1)(10^3-10)}{11^2}=8190.
$$

{\bf Example 4.2.} We use Theorem 1.2 to compute the number $N^{(n_1,n_2,n_3)}$ of
rational points on the following variety over $\mathbb{F}_7$:
\begin{align*}
{\left\{\begin{array}{rl}
x_1x_2^2x_3^2+x_1^2x_2x_3^5x_4^3x_5+x_1x_2^4x_3^3x_4^2x_5^4x_6x_7=1,\\
x_1^3x_2^2x_3^5+x_1x_2^3x_3^5x_4x_5^2+x_1x_2^2x_3x_4^4x_5^3x_6^2x_7^3=3,\\
x_1^2x_2^5x_3^2+x_1^2x_2^2x_3^4x_4^3x_5+x_1x_2^3x_3x_4^4x_5^3x_6^2x_7=5.\\
\end{array}\right.}
\end{align*}

Clearly, we have $n_1=3$, $n_2=5$, $n_3=7$, $r_1=1$, $r_2=2$, $r_3=3$, $m=3$,
$$E^{(1)}=\left( \begin{array}{*{20}c}
1&2&2\\
3&2&5\\
2&5&2\\
\end{array}\right),
$$
$$
E^{(2)}=\left( \begin{array}{*{20}c}
1&2&2&0&0\\
2&1&5&3&1\\
3&2&5&0&0\\
1&3&5&1&2\\
2&5&2&0&0\\
2&2&4&3&1\\
\end{array}\right) \  \text{and} \
E^{(3)}=\left( \begin{array}{*{20}c}
1&2&2&0&0&0&0\\
2&1&5&3&1&0&0\\
1&4&3&2&4&1&1\\
3&2&5&0&0&0&0\\
1&3&5&1&2&0&0\\
1&2&1&4&3&2&3\\
2&5&2&0&0&0&0\\
2&2&4&3&1&0&0\\
1&3&1&4&3&2&1\\
\end{array}\right).$$

We first calculate $N_{r_1}^{(n_1,n_2,n_3)}$.
Using elementary transformations, we obtain two unimodular matrices

$$U^{(1)}=\left( \begin{array}{*{20}c}
1&0&0\\
-2&0&1\\
11&-1&-4
\end{array}\right)
\ \text{and} \
V^{(1)}=\left( \begin{array}{*{20}c}
1&-2&-6\\
0&1&2\\
0&0&1
\end{array}\right)$$
such that
$$U^{(1)}E^{(1)}V^{(1)}={\rm SNF}
(E^{(1)})=\left( \begin{array}{*{20}c}
1&0&0\\
0&1&0\\
0&0&9\\
\end{array}\right).$$
Thus $d_1^{(1)}=1$, $d_2^{(1)}=1$, $d_3^{(1)}=9$ and $s_1=3$.
Clearly, the vectors $(u_{11},u_{21},u_{31})\in(\mathbb{F}^*_7)^3$ such that
\begin{align*}
{\left\{\begin{array}{rl}
u_{11}=1\\
u_{21}=3\\
u_{31}=5
\end{array}\right.}
\end{align*}
are
$$(u_{11},u_{21},u_{31})=(1,3,5).$$
Choose the primitive element 3 of $\mathbb{F}^*_7$. Then we have
$$
(h'_1, h'_2)^T=U^{(1)}(\text{ind}_31,\text{ind} _33,
\text{ind} _35)^T=U^{(1)}(6,1,5)^T\equiv(0,5,3)^T \pmod6.
$$
We deduce that the conditions (1.6) that $\gcd(6,9)|3$ hold.
It follows that $N_1=1$. So
\begin{align*}
N_{r_1}^{(n_1,n_2,n_3)}&=N_1q^{n_3-n_2}(q-1)^{n_1-s_1}(q^{n_2-n_1}-
(q-1)^{n_2-n_1})\prod\limits_{j=1}^{s_1}\gcd(q-1,d_j^{(1)})\\
&=7^2\times(7^2-6^2)\times3=1911.
\end{align*}

Consequently, we turn our attention to the computation of
$N_{r_2}^{(n_1,n_2,n_3)}$.
Using the elementary transformations, one gets that
$$U^{(2)}=\left( \begin{array}{*{20}c}
 1&0&0&0&0&0\\
 1&1&-1&0&0&0\\
 -2&-1&1&1&0&0\\
 8&6&-6&-4&0&1\\
 -1&5&-2&-1&2&-3\\
 7&-9&1&0&-5&9
\end{array}\right)
\ {\rm and} \
V^{(2)}=\left( \begin{array}{*{20}c}
1&-2&2&0&10\\
0&1&-2&-1&10\\
0&0&1&1&-15\\
0&0&0&0&1\\
0&0&0&-1&17
\end{array}\right)$$
such that
$$U^{(2)}E^{(2)}V^{(2)}
={\rm SNF}(E^{(2)})=\left( \begin{array}{*{20}c}
1&0&0&0&0\\
 0&1&0&0&0\\
 0&0&1&0&0\\
 0&0&0&1&0\\
 0&0&0&0&5\\
 0&0&0&0&0
\end{array}\right).$$
Thus $d_1^{(2)}=d_2^{(2)}=d_3^{(2)}
=d_4^{(2)}=1$, $d_5^{(2)}=5$ and $s_2=5$.
By Lemma 2.4, one has that the number of the vectors
$(u_{11},u_{12},u_{21},u_{22},u_{31},u_{32})\in(\mathbb{F}^*_7)^6$
such that
\begin{align*}
{\left\{\begin{array}{rl}
u_{11}+u_{12}=1\\
u_{21}+u_{22}=3\\
u_{31}+u_{32}=5
\end{array}\right.}
\end{align*}
is 125. Choose the primitive element 3 of $\mathbb{F}^*_7$.
The argument for calculating $N_1$ and using Matlab
we compute that $N_2=21$. Hence
\begin{align*}
N_{r_2}^{(n_1,n_2,n_3)}&=N_2q^{n_3-n_3}(q-1)^{n_2-s_2}(q^{n_3-n_2}-
(q-1)^{n_3-n_2})\prod\limits_{j=1}^{s_2}\gcd(q-1,d_j^{(2)})\\
&=21\times(7^2-6^2)=273.
\end{align*}

Let us now calculate $N_{r_3}^{(n_1,n_2,n_3)}$.
Using the elementary transformations, we obtain that
$$U^{(3)}=\left( \begin{array}{*{20}c}
1&0&0&0&0&0&0&0&0\\
-1&0&0&0&1&0&0&0&0\\
-2& -1&0&1&1&0&0&0&0\\
2&0&1&0&-1&0&-1&0&0\\
-7&-1&-2&1 &3&1&2&0&0\\
-13&-3&-3&2&5&0 &3&2&1\\
-15&1&0&2&3&0&4&-2&0\\
16&0&4&-2&-6&-1&-5&2&-1\\
7&-9&0&1&0&0&-5&9&0
\end{array}\right)$$
and
$$V^{(3)}=\left( \begin{array}{*{20}c}
1&-2&4&10&10&10&-20\\
0&1&-3&-7&-7&-7&15\\
0&0&1&2&2&2&-5\\
0&0&0&1&1&1&-2\\
0&0&0&0&0&0&1\\
0&0&0&0&-2&-1&0\\
0&0 &0&0&1&0&0
\end{array}\right)$$
such that

$$U^{(3)}E^{(3)}V^{(3)}
={\rm SNF}(E^{(3)})=\left( \begin{array}{*{20}c}
1&0&0&0&0&0&0\\
0&1&0&0&0&0&0\\
0&0&1&0&0&0&0\\
0&0&0&1&0&0&0\\
0&0&0&0&1&0&0\\
0&0&0&0&0&1&0\\
0&0&0&0&0&0&5\\
0&0&0&0&0&0&0\\
0&0&0&0&0&0&0\\
\end{array}\right).$$
Thus $d_1^{(3)}=d_2^{(3)}=d_3^{(3)}
=d_4^{(3)}=d_5^{(3)}=d_6^{(3)}=1$,
$d_7^{(3)}=5$ and $s_3=7$.
Lemma 2.4 gives that the number of the vectors
$(u_{11},u_{12},u_{13},u_{21},u_{22},u_{23},u_{31},u_{32},u_{33})
\in (\mathbb{F}^*_7)^9$ such that
\begin{align*}
{\left\{\begin{array}{rl}
u_{11}+u_{12}+u_{13}=1\\
u_{21}+u_{22}+u_{23}=3\\
u_{31}+u_{32}+u_{33}=5
\end{array}\right.}
\end{align*}
is equal to 29791. Choose the primitive element 3 of $\mathbb{F}^*_7$.
By the argument for calculating $N_1$ and using Matlab, we compute that
$N_3=823$. Thus one has
$$
N_{r_3}^{(n_1,n_2,n_3)}
=N_3(q-1)^{n_3-s_3}\prod\limits_{j=1}^{s_3}\gcd(q-1, d_j^{(3)})
=823.
$$
Finally, by Theorem 1.2, we have
$$N^{(n_1,n_2,n_3)}=\sum\limits_{i=1}^{3}N_{r_i}^{(n_1,n_2,n_3)}=1911+273+823=3007.$$


\end{document}